\documentclass{article}

\usepackage{fullpage}
\usepackage{epsfig}
\usepackage{graphicx}
\usepackage{amsmath}
\usepackage{amssymb}
\usepackage{color}
\usepackage{amsthm}
\usepackage{todonotes}
\usepackage{bbm}

\newcommand{\R}{\mathbb R}

\newcommand{\mB}{\mathcal B}

\newcommand{\mS}{\mathcal S}

\newcommand{\mA}{\mathcal A}
\newcommand{\mI}{\mathcal I}

\newcommand{\mH}{\mathcal H}

\newcommand{\tr}{\mathbf{tr}}

\newcommand{\diag}{\mathbf{diag}}
\newcommand{\find}{\mathbf{find}}

\newcommand{\proj}{\mathbf{proj}}
\newcommand{\evec}{\mathbf{evec}}

\newcommand{\range}{\mathbf{range}}
\newcommand{\vnull}{\mathbf{null}}

\newcommand{\st}{\mathbf{st}}

\newcommand{\mb}[1]{\mathbf{#1} }

\newcommand{\minimize}[1]{\underset{#1}{\textbf{minimize}}}
\newcommand{\mini}[1]{\underset{#1}{\textbf{min}}}

\newcommand{\argmax}[1]{\underset{#1}{\textbf{argmax}}}
\newcommand{\subjectto}{\mathbf{subject\; to}}

\newcommand{\mfind}{\mathbf{find}}
\newcommand{\rank}{\mathbf{rank}}

\newcommand{\bmat}{\left[\begin{matrix}}
\newcommand{\emat}{\end{matrix}\right]}

\usepackage{algorithm,algpseudocode}% http://ctan.org/pkg/{algorithms,algorithmx}

\algnewcommand{\Inputs}{%
	\State \textbf{Inputs:}
}
\algnewcommand{\Initialize}{%
	\State \textbf{Initialize:}
}
\algnewcommand{\Outputs}{%
	\State \textbf{Outputs:}
}

\begin{document}

%%%%%%%%% TITLE
\title{Approximate methods for phase retrieval via gauge duality}
\author{Ron Estrin, Yifan Sun, Halyun Jeong, Michael Friedlander}
\maketitle

\begin{abstract}
We consider the problem of finding a low rank symmetric matrix satisfying a system of linear equations, as appears in phase retrieval. In particular, we solve the gauge dual formulation, but use a fast approximation of the spectral computations to achieve a noisy solution estimate. This estimate is then used as the initialization of an alternating gradient descent scheme over a nonconvex rank-1 matrix factorization formulation. Numerical results on small problems show consistent recovery, with very low computational cost.
\end{abstract}

\section{Introduction}
Consider the problem of finding a low rank symmetric matrix satisfying a system of linear equations
\begin{equation}
\begin{array}{ll}
\mfind & X \\
\subjectto & a_i^TX a_i = b_i,\quad i = 1,\hdots, m\\
& \rank(X) \leq r
\end{array}
\label{e-main}
\end{equation}
Problems of form \eqref{e-main} appear in applications like 
imaging \cite{walther1963question} and x-ray crystallography \cite{dierolf2010ptychographic}, and finding $x$ is in general NP-hard \cite{vavasis2010complexity}. 
Convex relaxations of \eqref{e-main} appears by omitting the rank constrant, and can often lead to a close approximation of $x$ (see  \cite{ahmed2014blind,candes2013phaselift}).
 
We consider two convex relaxations of \eqref{e-main}, both of the form
\begin{equation}
\begin{array}{ll}
\minimize{X}& \kappa(X)\\
\subjectto & a_i^TX a_i = b_i,\quad i = 1,\hdots, m
\end{array}
\label{eq:main-primal}
\end{equation}
where $\kappa(X)$ is a convex gauge function that promotes low-rank structure in $X$. Specifically, we consider two choices:
\begin{enumerate}
\item $\kappa(X) = \|X\|_*$ the nuclear norm, e.g. the sum of the singular values of $X$, and 
\item $\kappa(X) = \tr(X) + \delta_+(X)$ where
\[
\delta_+(X) = 
\begin{cases}
0 & X\succeq 0 \\
+\infty & \text{ else.}
\end{cases}
\]
\end{enumerate}
Using either gauge, we see that \eqref{eq:main-primal} is a semidefinite optimization problem.

\paragraph{Application 1: Phase retrieval}
Problems of this form appears in image processing as a convexification of the phase retrieval problem 
\begin{equation}
\begin{array}{ll}
\find& x\\
\subjectto & |a_i^Tx|^2 = b_i,\quad i = 1,\hdots, m 
\end{array}
\label{eq:phaseretrieval}
\end{equation}
where $a_i$ may be complex or real-valued measurements, and $b_i$ are the squared magnitude readings. 
 (See PhaseLift, \cite{candes2015phase}.) In particular, it has been shown \cite{chen2001atomic} that low-rank estimates of the SDP \eqref{eq:main-primal} can recover the exact source vector $x$ in both the noisy and exact measurement regime, for large enough $m$ and incoherent enough $a_i$. 
 It has been previously observed (and numerically verified here) that the positive semidefinite formulation ($\kappa(X) = \tr(X) + \delta_+(X)$) provides better recovery results (recovering successfuly for smaller $m$) than the nuclear norm formulation. However, we will see there are numerical advantages of using the nuclear norm formulation, and thus we consider both. 
 
\paragraph{Application 2: Linear diagonal constrained}
Many combinatorial problems can be relaxed to semidefinite programs with diagonal constraints. For example,  the MAX-CUT problem can be written as 
 \begin{equation}
\begin{array}{ll}
\minimize{W}& \langle C,W\rangle\\
\subjectto & \diag(W) = b \\
& W\succeq 0
\end{array}
\label{eq:diag-cons}
\end{equation}
for $b = \mb 1$ and  $C$ is a matrix related to the graph edge weights (e.g. Laplacian).
Considering the a more generalized family of linear diagonally constrained problems, we see that replacing $C$ with $(C+C^T)/2 + \diag(v)$ does not alter the problem, for any $v$, since the diagonal of $W$ is fixed. Therefore we can assume without loss of generality that $C$ is symmetric positive definite with Cholesky factorization $C = LL^T$. Then \eqref{eq:diag-cons} is equivalent to \eqref{eq:main-primal} where $a_i$ is the $i$th column of $L^{-1}$.

\section{Problem statement} 
\subsection{Gauge duality}
In general, the gauge primal and dual problem pair\cite{freund1987dual,friedlander2014gauge,aravkin2017foundations,friedlander2016low} can be written as
\begin{equation}
\begin{array}[t]{ll}
\mini{X}& \displaystyle \kappa(X)\\
\st & \mA(X) = b,\quad i = 1,\hdots, m
\end{array}
\qquad
\begin{array}[t]{ll}
\mini{y}& \displaystyle \kappa^\circ(\mA^*(y))\\
\st & \langle y,b\rangle = 1
\end{array}
\label{eq:gauge-pair}
\end{equation}
 where we use the shorthand
\[
\mA(X)_i := a_i^TXa_i,\; i = 1,\hdots, m, \qquad \mA^*(y):=\sum_{i=1}^m y_i a_ia_i^T
\]
for the linear operator and adjoint.
Here, $\kappa^\circ$ is the polar gauge of $\kappa$, defined as
\[
\kappa^\circ(Z) = \inf \{\mu : \langle X, Z\rangle \leq \mu \kappa(X)\; \forall X. \} 
\]
In particular, it is shown that if the feasible domain of both primal  and dual \eqref{eq:gauge-pair} have nontrivial relative interior, then at optimality, the eigenspace of the primal matrix variable $X^*$ and transformed dual variable $Z^* = \mA^*(y*)$ are closely related, and can often be recovered easily.

%\blue{See Corollary 5.2 in \cite{friedlander2014gauge} for more details. I guess the similar properties for the Lagrange dual are true; Slater's condition $\Rightarrow$ strong duality $\Rightarrow$ strictly complementary $\Rightarrow$ $\range(X^*) = \textbf{null}(Z^*)$. See \cite{ding2019optimal} based on this relation (Actually it is general since they used approximate complementary).}

\paragraph{Nuclear norm}
When $\kappa(X)$ corresponds to a norm, then $\kappa^\circ(X)$ is the dual norm. Therefore 
\[
\kappa(X) = \|X\|_* \iff \kappa^\circ (Z) = \|Z\|_2
\]
the spectral norm of $Z$. Note that neither $X$ nor $Z$ are constrained to be positive semidefinite. 
At optimality, the singular vectors of the primal matrix variable $X^*$ and transformed dual variable $Z^* = \mA^*(y*)$ correspond closely;  if $X^*$ has rank $r$, then
\[
X^* = \sum_{i=1}^r \sigma^P_i v_iv_i^T, \qquad Z^* = \sigma^D_{\max}\sum_{i=1}^r v_iv_i^T + \sum_{i=r+1}^n \sigma^D_i v_iv_i^T.
\]
Here, $v_1,\hdots, v_n$ are the singular vectors of $X^*$ and $Z^*$, and $\sigma^P_i$ and $\sigma^D_{\max}$ correspond to the primal singular values and maximum dual singular values, respectively. 
Note that the singular vectors of the primal and dual variables are the \emph{same}, so the range of $X^*$ can be recovered from either primal or dual optimal solutions.

%\blue{True. Sometimes we don't have exact dual optimal solution but still can recover $\hat{X}$ close to $X^*$. }

\paragraph{Symmetric PSD}
In the second case, 
\[
\kappa(X) = \tr(X) + \delta_+(X) \iff \kappa^\circ (Z) = \max\{0,\lambda_{\max}(Z)\}.
\]
Through gauge duality, $X^*$ and $Z^* = \mA^*(y^*)$ 
 have a \emph{simultaneous eigendecomposition}; that is, if $X^*$ has rank $r$, then
\[
X^* = \sum_{i=1}^r \lambda^P_i u_iu_i^T, \qquad Z^* = \lambda^D_{\max}\sum_{i=1}^r u_iu_i^T + \sum_{i=r+1}^n \lambda_i u_iu_i^T.
\]
Here, $u_i$ are the eigenvectors of both $X^*$ and $Z^*$, and $\lambda^P_i$ and $\lambda^D_{\max}$ correspond to the primal eigenvalues and maximal dual eigenvalue, respectively.

Additionally, strong  gauge duality enforces $1 = \kappa(X)\kappa^\circ(Z)$ at optimality. Assuming the primal of \eqref{eq:gauge-pair} is feasible, $\kappa(X) < +\infty$ which forces $\kappa^\circ(Z^*) > 0$. Therefore, we can simplify the dual objective function to
\[
\kappa^\circ (Z) = \lambda_{\max}(Z)
\]
over $Z$ where  $\lambda_{\max}(Z) > 0$.
\paragraph{Unconstrained formulation}

We now rewrite the dual of \eqref{eq:gauge-pair}  in an unconstrained formulation
 \begin{equation}
\minimize{z}\quad \displaystyle \kappa^\circ(\mA^*( Bz+\bar z ))
\label{eq:main-dual-uncons}
\end{equation}
using a change of variables $y = Bz+\bar y$ for any $B$ where $\range(B) = \vnull(b^T)$ and $\bar y$ such that $\langle \bar y,b\rangle = 1$. In this case, for any $z$, 
$\langle y,b\rangle = 1$.

\section{Methods}

\subsection{General overview}
We consider three methods, described in ``vanilla" form below.
\begin{enumerate}
\item \emph{Projected gradient descent}  on the constrained gauge dual of \eqref{eq:gauge-pair}
\begin{equation}
		y^{(k)}  = \proj_{\mH}(y^{(k-1)} - t\nabla_y \kappa^\circ (y^{(k-1)}))
\end{equation}
where $\mH = \{y : \langle y,b\rangle = 1\}$ is the constraint set. The Euclidean projection on this set can be done efficiently via
\[
		\proj_{\mH}(s) =  \left(I-\frac{1}{b^Tb} bb^T\right) s + \frac{1}{b^Tb} b.
\]

\item \emph{Gradient descent} on the unconstrained gauge dual \eqref{eq:main-dual-uncons}
\begin{equation}
		z^{(k)}  = z^{(k-1)} - tB^T g
		\label{eq:unconstrained-descent-step}
\end{equation}
where
\[
g = \nabla_y \kappa^\circ (y) \quad \text{ at }\quad y = Bz^{(k-1)} + \bar y.
\]
An obvious choice for $B\in \R^{m-1,m}$ computed from a full QR of $b^T$ where $B^TB = I$ and $Bb = 0$. 
\item \emph{Coordinate descent} (Alg. \ref{a:cd}) on the unconstrained gauge dual \eqref{eq:main-dual-uncons}
\begin{equation}
		z^{(k)}  = z^{(k-1)} - tB^T g
\end{equation}
where
\[
(g_z)_j =
		\begin{cases}
		\frac{\partial \kappa^\circ(y)}{\partial z_j} & j\in \widehat \mI^{(k)}\\
		0 & \text{ else.}
		\end{cases}
		 \quad \text{ at }\quad y = Bz^{(k-1)} + \bar y.
\]
and $\widehat \mB^{(k)} = \{j : i\in \mB^{(k)}, B_{ij} \neq 0\}$.
Here, in order to maintain scalability, we want $|\widehat \mB^{(k)}|$ to be small whenever $|\mB^{(k)}|$ to be small, e.g. $B$ should be sparse. Note $B$ does not have to be orthogonal  for the unconstrained formulation \eqref{eq:main-dual-uncons} to be equivalent to the dual constrained formulation \eqref{eq:gauge-pair}. However, we have found much better results when $B$ is as close to $I$ as possible, so we pick 
\begin{equation}
B = \bmat I_i  & \tilde b_{\mathrm{top}} & 0 \\0  & \tilde b_{\mathrm{bot}} & I_{m-i} \emat, 
\label{eq:sparseB}
\end{equation}
where
\[
 \tilde b_{\mathrm{top}} = \left(\frac{b_1}{b_i},\hdots,  \frac{b_{i-1}}{b_i}\right), \qquad \tilde b_{\mathrm{bot}} = \left(\frac{b_{i+1}}{b_i},\hdots, \frac{b_m}{b_i}\right)
\]
and $i = \argmax{i}\; b_i$. 

%\blue{TODO: Estimate the condition number of $B$ on $\textbf{null}(b^\perp)$ when $B$ is close to $I$. This may tell us why we have better result in this case.}

\end{enumerate}
The main focus of this work is to exploit structural properties of the linear operator $\mA$, and offer several ``enhancements" that significantly improve the scalability of these methods. In particular, we will focus on 

\begin{enumerate}
\item approximations for building the gradients $\nabla_y \kappa^\circ(y)$  (or $s\in \partial_y \kappa^*(y)$ in cases where the largest eigenvalue of the formed matrix is not simple)
\item picking $B$ in the unconstrained formulation so that multiplying by $B$ and $B^T$ is efficient, and

\item estimating the partial coordinates in the coordinate descent method. 
This is our primary contribution is the third improvement, which theoretically can avoid doing any spectral computations (\texttt{svds} or \texttt{eigs}) and is limited only to small matrix products and a tiny QR computation, and maintains only low-rank approximations of all matrices. 
\footnote{In practice, in order to reach the global solution, sometimes the matrix estimates deteriorate, and need to be ``refreshed", so a full \texttt{eigs} is run. However, in phase retrieval we often don't need this much precision.}

\end{enumerate}

\subsection{Gradients of dual objective}
Let us first consider 
\[
f(y) = \kappa^\circ(\mA^*(y)) = \max\{\lambda_{\max}(\mA^*(y)),0\}.
\]
 Then computing $\nabla f(y)$ requires three steps.
\begin{enumerate}
\item Form the dual matrix variable 
\[
W = \mA^*(y)
\]

\item Find $u$ where $Wu = \lambda_{\max}u$. Since this step is important, we will denote this operation $u = \evec_{\max}(W)$.

\item
The gradient is now
\[
\nabla f(y) = 
\bmat 
(a_1^T u)^2\\
\vdots\\
(a_m^T u)^2.
\emat
\]
This step is comparatively cheap, so we will not discuss it.
\end{enumerate}
Now consider $m$ and $n$ both large, with $m > n$. 

\paragraph{Building $W$} When $m$ is large, this step is both computationally expensive ($O(n^2 m)$) and memory inefficient if $n$ is large. 
Note that if $y > 0$, there is a simple computational shortcut is to form 
\[
U = \bmat \sqrt{y_1} a_1,\hdots , \sqrt{y_m}a_m\emat, \qquad W = UU^T.
\]
However, in general the intermediate $y^{(k)}$ and final $y^*$ are not nonnegative. 
We therefore try to estimate this quantity using
\begin{equation}
\widehat W = \sum_{i\in \mB} y_i a_ia_i^T
\label{eq:approx-matrix}
\end{equation}
and $\mB$ is some sample subset of $\{1,...,m\}$. 
In particular, we investigate three regimes
\begin{enumerate}
\item $\mB = \{1,...,m\}$ exact gradient computation
\item $\mB = \{i : y_i \geq 0 \}$ for a PSD estimate of $\mA^*(y)$ 
\item $j\in \mB$ is randomly picked from the set $\{i : y_i > 0\}$ with probability $y_i / \sum_{y_j \geq 0}y_j$. 
\end{enumerate}

\begin{figure}
	\begin{center}
\includegraphics[width=3in]{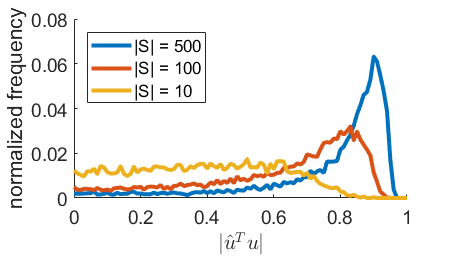}
\caption{\textbf{Error in subsampling.} In this experiment, we consider $Z = \sum_{i=1}^m y_ia_ia_i^T$ where $m = 1000$ and $a_i\in \R^n$ and $y_i\in \R$ are i.i.d. Gaussian randomly generated vectors and scalars. Sorting the weights $y_{i_1} \geq y_{i_2} \geq \cdots$, we select  $\mS = \{i_1,...,i_{|S|} \}$ and $\hat Z = \sum_{i\in \mS} y_ia_ia_i^T$, and compare the alignment between $u = \evec_{\max}(Z)$ and $\hat u = \evec_{\max}(\hat Z)$.}
\label{fig:subsamp-grad}
\end{center}
\end{figure}

\paragraph{Solving the eigenvalue problem} 
This step can be solved fairly quickly using a fast eigenvalue solver (such as a power method). However, a key issue when $W$ is indefinite is that the largest magnitude eigenvalue may be much larger than the largest algebraic eigenvalue. 
This is another key motivation behind the second choice of $\mB$, to work with a PSD estimate $\widehat W$. 
In practice, we do not observe performance degradation with this estimation, and in fact observe considerable speedup in the convergence of \texttt{eigs}. 
(See also Figure \ref{fig:subsamp-grad}.)

\paragraph{Extension to nuclear norm}
The above discussion mostly also holds for $f(y) = \kappa^\circ(\mA^*(y))  = \|\mA^*(y)\|_*$, but replacing eigenvalue computations with singular value computations, and sampling based on $|y_i|$ rather than $y_i$.  A key subtle advantage of using the nuclear norm is that since singular values are always nonnegative, we do not need to worry as much about the convergence of \texttt{svds}. Note that though nuclear norm minimization is generally used for nonsymmetric matrix variables, here because $\mA^*(y)$, we are still only considering symmetric matrix variables. 
(The distinction is that we are now running \texttt{eigs(W,1,`lm')} where previously we ran \texttt{eigs(W,1,`la')}.)

\subsection{Coordinate descent}
For applications where both $m$ and $n$ are large, we further parametrize $W$ with a low-rank approximation and use a block coordinate update at each iteration. 
This method is inspired by the following observation: if $z$ and $\hat z$ differ by at most $k$ elements, then for $B$ constructed as in \eqref{eq:sparseB},
\[
y = Bz + \bar y \text{ and } \hat y = B\hat z + \bar y
\]
differ by at most $2k$ elements, and 
\[
W = \mA^*(y) \text{ and } \widehat W = \mA^*(\hat y)
\]
differ by at most a term of rank $2k$. Now assume that at iteration $k$, we maintained a rank-$r$ approximation of $W^{(k)}$ as
\begin{eqnarray*}
U^{(k)}D^{(k)} (U^{(k)})^T &=&  U^{(k-1)}D^{(k-1)} (U^{(k-1)})^T + A_{\mI^{(k)}} \diag(y_{\mI^{(k)}}) (A_{\mI^{(k)}})^T\\
&=&  \bmat U^{(k-1)} &  A_{\mI^{(k)}}\emat \bmat D^{(k-1)}  & 0 \\ 0 & \diag(y_{\mI^{(k)}}) \emat \bmat (U^{(k-1)})^T  \\ (A_{\mI^{(k)}})^T\emat 
\end{eqnarray*}
where the columns of $A_{\mI^{(k)}}$ are $a_i$ for $i\in \mI^{(k)}$. Packing $\widetilde D = \bmat D^{(k-1)}  & 0 \\ 0 & \diag(y_{\mI^{(k)}}) \emat$, taking a QR factorization of $[U^{(k-1)}, A_{\mI^{(k)}}] = QR$, we have
\begin{eqnarray*}
U^{(k)}D^{(k)} (U^{(k)})^T &=&  
  QR \widetilde D  R^TQ^T
\end{eqnarray*}
where $R\widetilde D R^T$ is $r+2k \times r+2k$, and in general $r+2k \ll m$. Taking a ``tiny eig" of this matrix
\[
R\widetilde D R^T = \widetilde U  D^{(k)} \widetilde U^T
\]
gives the new rank-$r+2k$ factorization of $W^{(k)}$ with $U^{(k)} = Q\widetilde U$. In the algorithm, we then further prune $D^{(k)}$ and $U^{(k)}$ to its best rank-$r$ approximation.

\paragraph{Picking the coordinates}
At each iteration, the coordinates $\mI^{(k)}$ can be picked uniformly without replacement from $\{1,\hdots, m\}$, or according to a greedy method. In particular, the Gauss-Southwell ``flavor" of coordinate descent algorithms picks 
\[
i = \argmax{i}\; |(\nabla \kappa^\circ (y))_i|.
\]
However, note that just making this call requires computing a full gradient, which we never want to do. Therefore we approximate this operation by sampling at each iteration $i$ according to a \emph{weighted} uniform distribution, with weights $\max\{y_i,0\}$ when $\kappa = \tr + \delta_+$ and $|y_i|$ when $\kappa = \|\cdot\|_*$.

\begin{algorithm}
	\caption{Maintaining low rank steps}
	\begin{algorithmic}[1]
		\Inputs{rank parameter $r$}
		\Outputs {$z^{(k)}$}
		
		\State $\bar z = e_{i_{\max}}$ where $i_{\max} = \argmax{i}\;|b_i|$.
		
		\State $z^{(0)} = 0$, $y^{(0)} := \bar z$, $W^{(0)} := \mA^*(y^{(0)})$
		\State Compute top-$r$ eigenvalue decomposition
		\[
		UDU^T = \proj_{\rank = r}(W^{(0)})
		\]
		and the diagonal of $D$ is decreasing in order.
		\item Set $k = 0$.
		\For{$k = 1,...$}
			\State Sample $\widehat\mI^{(k)}$ containing $L$ elements without replacement from $\{1,\hdots, m\}$ and update
			\[
			\mI^{(k)} = \{i : B_{ij} \neq 0,\; i\in \widehat \mI^{(k)}\}.
			\]
			
			\State Compute partial gradients of $\kappa^\circ(\mA^*(y))$  with respect to $y$ and $z$, with  $u = U[:,1]$
        		\[
        		g_i  =  (a_i^T u)^2, \; i\in \mI^{(k)},\qquad
			\hat g_j = 
			\begin{cases}
			\displaystyle \sum_{i:B_{ij}\neq 0} B_{ij} g_i,& j\in \widehat\mI^{(k)} \\
			0 & \mathrm{ else.}
			\end{cases}
			\]

			\State Update $z^{(k)}$
			\[
			z^{(k)} := z^{(k-1)} - t \hat g.
			\]
			\State Update the rank $r$ approximation of $W^{(k)}$
			\[
			\widetilde D = \bmat D^{(k)} & 0 \\ 0 &  \diag(\Delta y) \emat ,\qquad 
			QR = \texttt{qr}(\bmat V &   A_{\mI^{(k)}}^T\emat,0),\qquad 
        \Delta y =  - t B^T\hat g 
			\]
			\[
			[\widetilde U,\widehat D] = \texttt{eig}(R*\widetilde D*R'),\qquad 
			\widehat U = Q \widetilde U
			\]
			\State Prune to rank $r$
			\[
			U^{(k)} D^{(k)} (U^{(k)})^T = \proj_{\rank=r}(\widehat U \widehat D {\widehat U}^T)
			\]
		\EndFor

	\end{algorithmic}
\label{a:cd}
\end{algorithm}

\section{Numerical results}

\subsection{Musical note}
\begin{figure}
	\begin{center}
	\fbox{\includegraphics[width=1in]{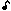}}
\caption{\textbf{Musical note.}}
\label{fig:musicalnote}
\end{center}
\end{figure}

We begin by considering a small, simple problem of recovering an 11 x 11 black and white image (Figure \ref{fig:musicalnote}). This problem is small enough to be solved globally using an interior point method, which can serve as a baseline. 

\begin{figure}
\centering
\includegraphics[width=0.7\linewidth]{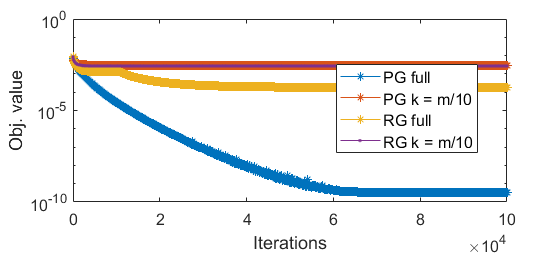}
\includegraphics[width=0.7\linewidth]{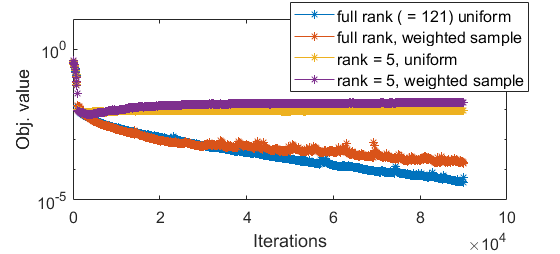}
\caption{\textbf{Trajectory.} Oversampled musical note example, with $m = 1000$, $n = 121$. \textbf{Top:} Gradient methods, comparing full gradients vs $m/10$ subsampled gradients. PG = projected gradient, RG = reduced gradient.
	\textbf{Bottom:} Coordinate methods, with block size 100. 
}
\label{fig:traj}
\end{figure}

\begin{figure}
\centering
\fbox{
\includegraphics[width=0.7\linewidth]{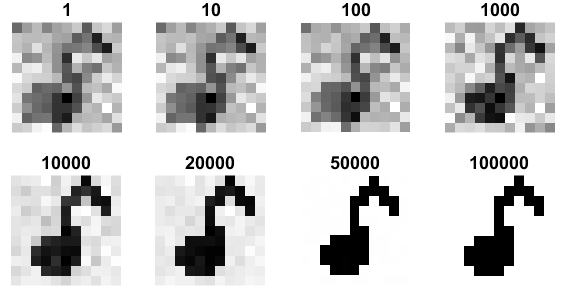}}
	\fbox{\includegraphics[width=0.7\linewidth]{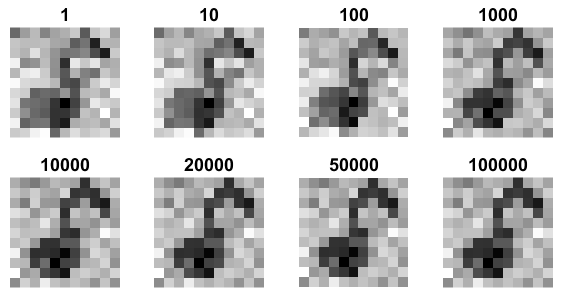}
	}
\caption{\textbf{Projected gradient descent.} Recovered image using $u^{(t)} = \evec_{\max}(Z^{(t)})$.Title gives iteration count $t$. \textbf{Top two rows:} full gradients used. \textbf{Bottom two rows:} $m/10$ gradients used.}
\label{fig:projgrad-images}
\end{figure}

\begin{figure}
\centering
\fbox{\includegraphics[width=0.7\linewidth]{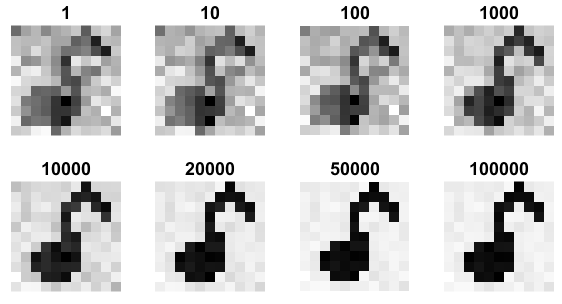}}
\fbox{\includegraphics[width=0.7\linewidth]{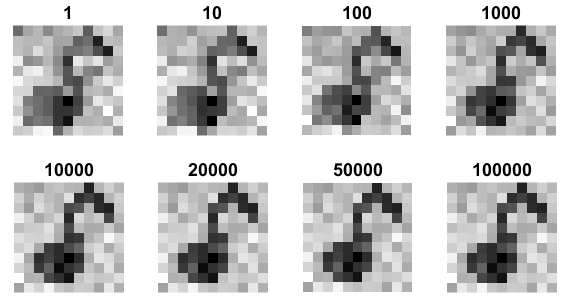}}
\caption{\textbf{Reduced  gradient descent.} Recovered image using $u^{(t)} = \evec_{\max}(Z^{(t)})$.Title gives iteration count $t$. \textbf{Top two rows:} full gradients used. \textbf{Bottom two rows:} $m/10$ gradients used.}
\label{fig:uncgrad-images}
\end{figure}

\begin{figure}
\centering
\fbox{\includegraphics[width=0.7\linewidth]{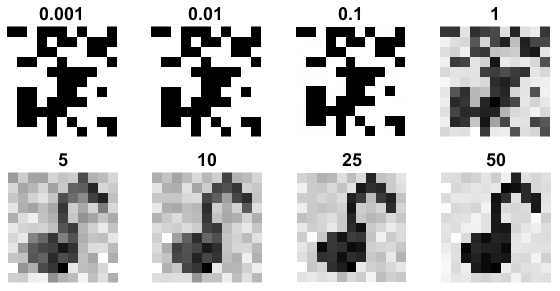}}
\fbox{\includegraphics[width=0.7\linewidth]{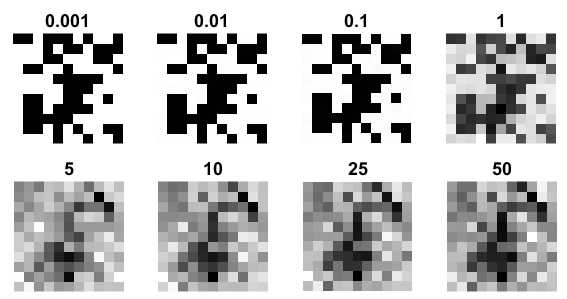}}
\caption{\textbf{Coordinate descent.} Recovered image using $u^{(t)} = \evec_{\max}(Z^{(t)})$. Title gives iteration count $t$. In all cases, block size = 100. \textbf{Top two rows:} full rank $Z^{(t)}$ stored. \textbf{Bottom two rows:} rank 5 approximation stored.}
\label{fig:coord-images}
\end{figure}

\paragraph{Fast methods do not give high enough fidelity solutions.} To evaluate how ``fast"  our methods work, we pick a fairly easy problem, with $m = 1000$ samples $a_i$ sampled uniformly without replacement from a Hadamard matrix. 
\begin{itemize}
	\item Figure \ref{fig:traj} shows the trajectory of the dual objective error for the first-order and coordinate methods on this problem. We can see that when full gradients and full rank methods are used, the global solution can be found, but the number of iteratious is onerous, especially since this is a very small example. When partial gradients and low rank approximations are used, the global solution is not found.

	\item Intermediate recovered images of the oversampled problem are given in Figure \ref{fig:projgrad-images} (projected gradient), \ref{fig:uncgrad-images} (reduced gradient), and \ref{fig:coord-images} (coordinate descent). Again, we notice that while non-approximated methods can recover the true solution (after many iterations), their approximated versions do not reach high fidelity solutions.

\end{itemize}

\paragraph{Fast methods approach good approximate solutions quickly.} One thing we do observe from the previous batch of experiments is that our fast methods are able to reach good approximate solutions almost immediately, suggesting they may provide good initializations to simpler nonconvex methods.

\paragraph{Nonconvex matrix completion}
Specifically, a common approach to phase retrieval is to solve the following nonconvex rank-1 matrix completion problem
\begin{equation}
\minimize{u} \quad \sum_{i=1}^m ((a_i^Tu)^2-b_i)^2
\label{eq:nonconvex-matrix-completion}
\end{equation}
Given some initial point $u^{(0)}$, we can minimize \eqref{eq:nonconvex-matrix-completion} by iteratively using gradient steps
\[
u^{(k+1)} = u^{(k)} -\alpha^{(k)} \sum_{i=1}^m ((a_i^Tu)^2-b_i) \cdot a_i
\]
where $\alpha^{(k)}$ is some decaying step size.
This type of approach is often favored in practice because of its low per-iteration complexity ($O(mn)$) and storage $O(n)$ cost, and is often observed to recover very clean images. A disadvantage of this approach is that the quality of the solution is very sensitive to the choice of initialization. As an example, the Wirtinger flow algorithm of \cite{candes2015phase} recovers images using the initialization 
\[
u^{(0)} = \evec_{\max}\left(\sum_{i=1}^m b_i a_ia_i^T\right).
\]

\begin{figure}
\begin{center}
\includegraphics[width=5.5in]{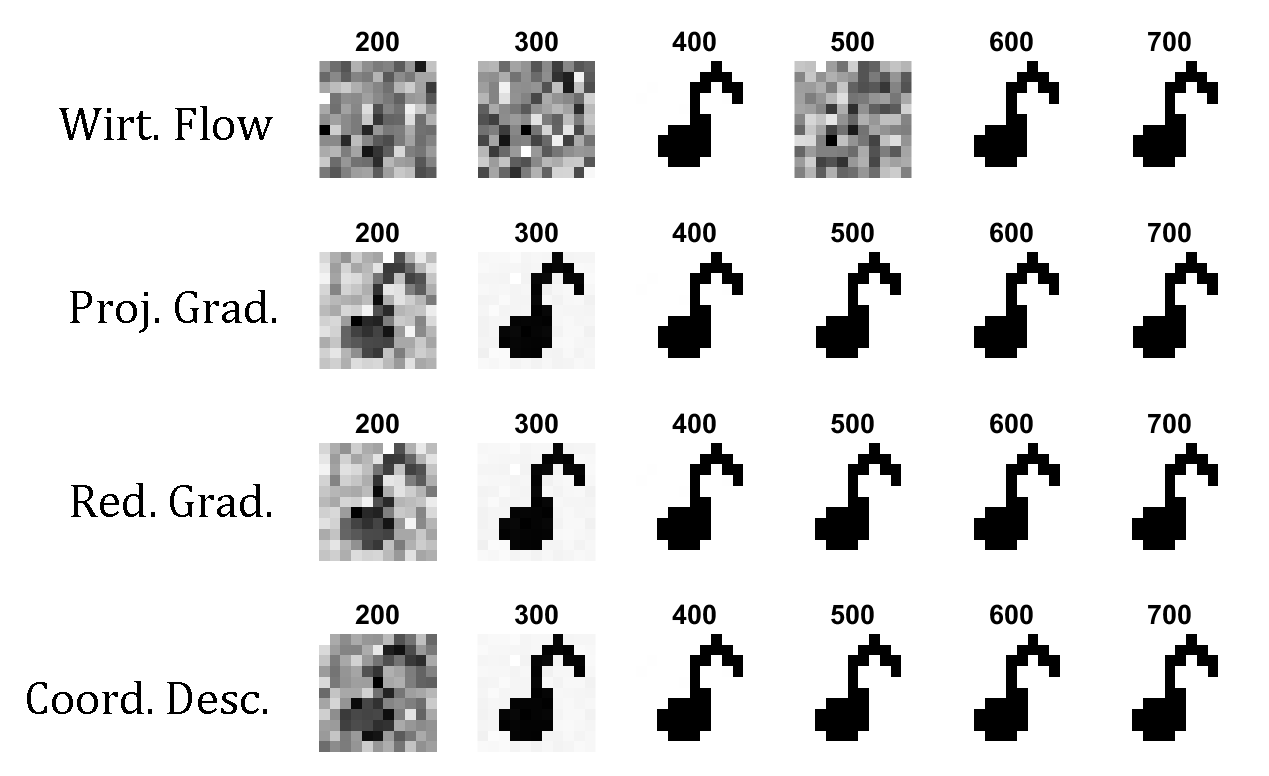}
\caption{\textbf{Solutions to nonconvex matrix completion.} Recovered images using a variety of initializations. Title = \# samples.}
\label{fig:musicnote-init}
\end{center}
\end{figure}

\begin{figure}
	\begin{center}
\includegraphics[width=5.5in]{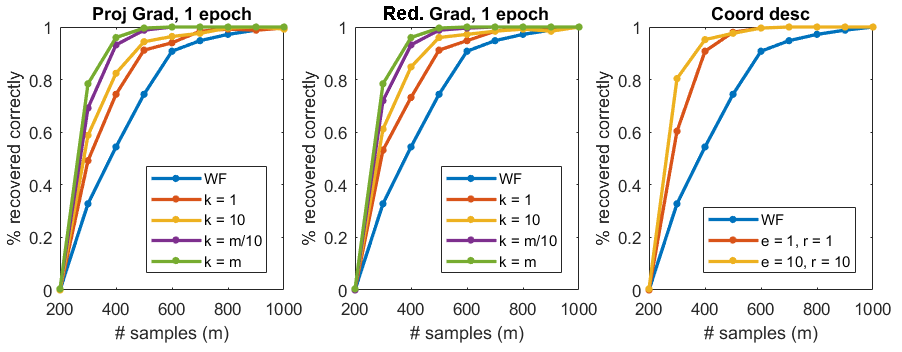}
		\caption{\textbf{Observed recovery rate.} k = \# samples in gradient. e = epochs, r = rank. WF = Wirtinger Flow initialization.}
		\label{fig:musicnote-recovery}
	\end{center}
\end{figure}

We propose to recover images using the nonconvex matrix completion method using, as initialization, an approximate solution from a few iterations of our faster methods
\[
u^{(0)} = \evec_{\max}\left(\sum_{i=1}^m y^{(K)}_i a_ia_i^T\right)
\]
Figure \ref{fig:musicnote-init} (random instance) and \ref{fig:musicnote-recovery} (averaged over 250 trials) illustrate the competitive advantage of using a solution to the fast method as an approximate solution as a higher fidelity initialization of the matrix completion problem.

\begin{figure}
\begin{center}
\includegraphics[width=6in]{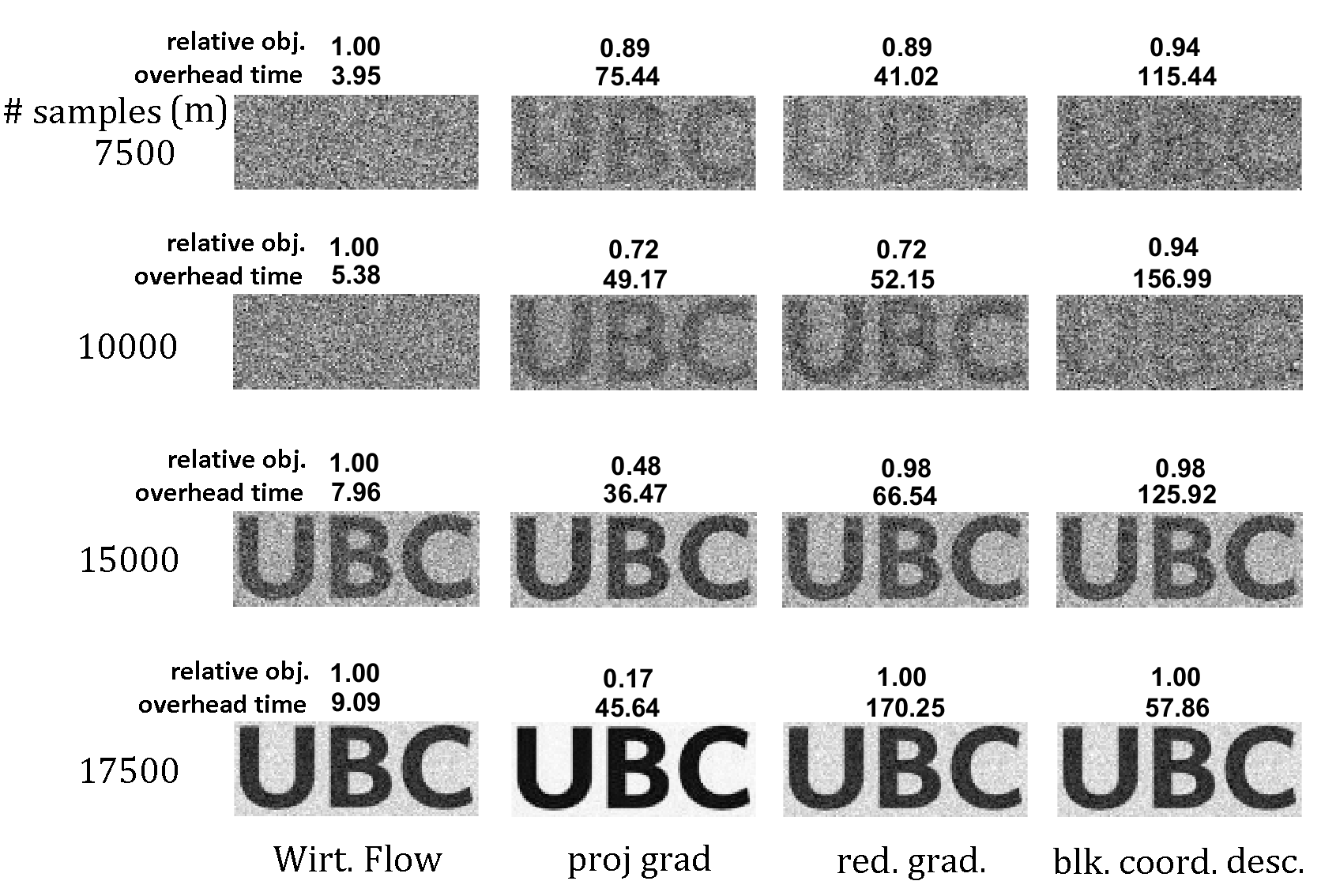}
\caption{\textbf{UBC Logo.} $n = 5220$. Relative objective = $f* / f*_{\mathrm{WF}}$. Overhead time refers to total time used to compute the initial point. The average runtime of the nonconvex matrix completion is about 30 seconds. All hyperparameters (number of iterations, step size decay scheme) are tuned to make each example as efficient and high quality as possible. }
\label{fig:ubclogo}
\end{center}
\end{figure}

\begin{figure}
	\begin{center}
\includegraphics[width=6in]{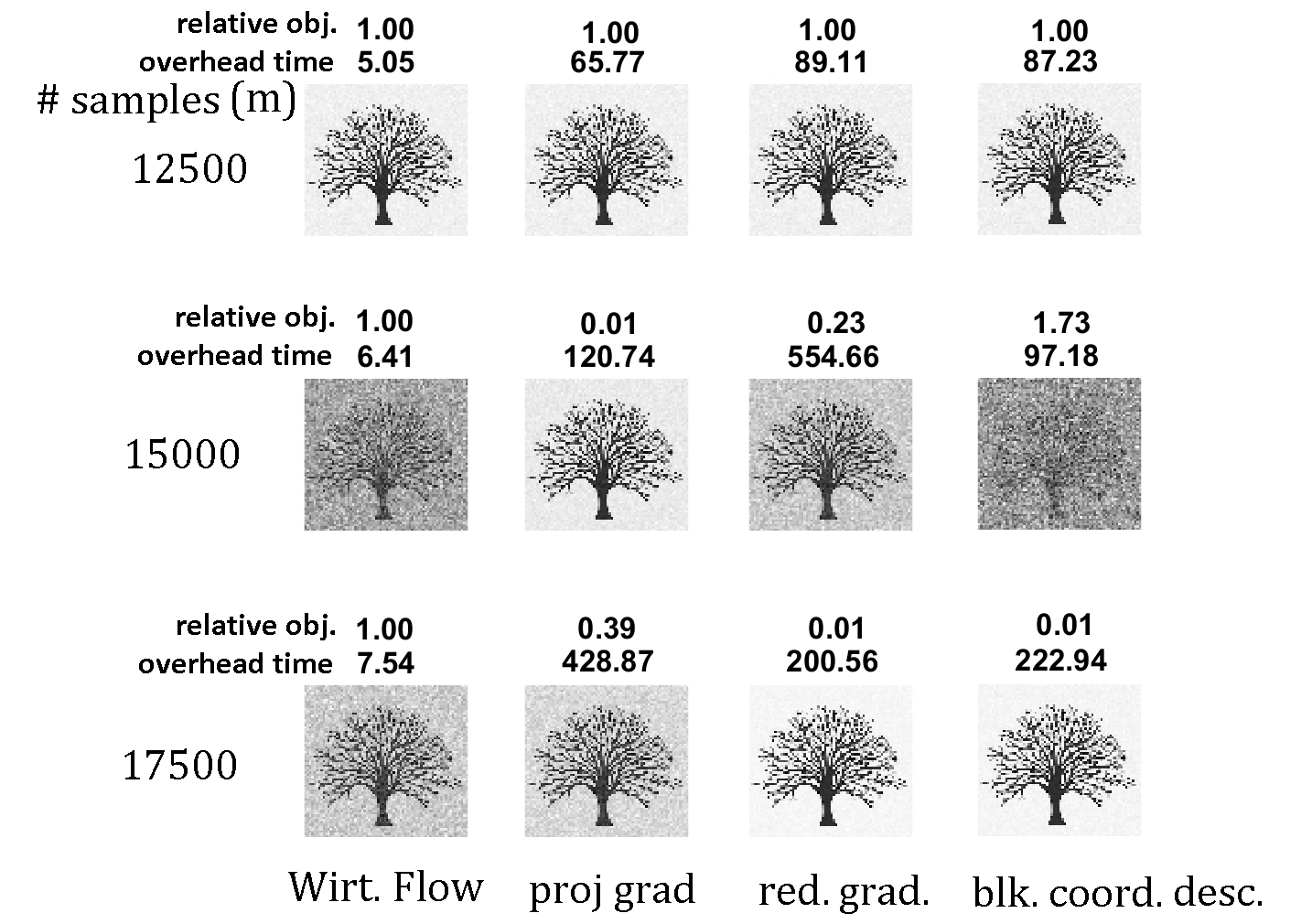}
		\caption{\textbf{Tree.} $n = 4824$. Relative objective = $f* / f*_{\mathrm{WF}}$. Overhead time refers to total time used to compute the initial point. The average runtime of the nonconvex matrix completion is about 5 minutes. All hyperparameters (number of iterations, step size decay scheme) are tuned to make each example as efficient and high quality as possible.}
		\label{fig:tree}
	\end{center}
\end{figure}

\paragraph{Slightly larger numerical results.} 
Figure \ref{fig:ubclogo} and \ref{fig:tree} repeat the experiment on slightly larger images, with different structural properties.

\section{Further directions}
This document represents a quick set of experiments on a simple idea for reducing the computational complexity of the SDP relaxation of the phase retrieval problem. There are several interesting directions for extension.

\paragraph{Gradient sampling}
Currently, our gradient sampling approach is to sample each weight according to its positive contribution, normalized, with no other transformations. A more generalized class of sampling weights is to use softmax smoothing, where 
\[
\mathrm{Pr}(j) = \frac{\exp(\frac{y_j}{\mu})}{\sum_k\exp(\frac{y_k}{\mu})}
\]
and for a specific choice of $\mu$, reduces to our sampling scheme. 
This kind of sampling  can be modeled using a Gumbel random variable, for example. 

\paragraph{Better choices of $B$} In the unconstrained dual formulation \eqref{eq:unconstrained-descent-step}, there is a tradeoff in the choice of $B$ as incredibly sparse (improving its per-iteration complexity) and perfectly conditioned (ideally, orthogonal, and therefore dense). Further exploration here can be made to optimize this tradeoff.

\paragraph{Scalability}
Thus far, we have viewed the most successes with very small images. With larger images, it is not clear how the approximation error scales, and if it is still close enough to ensure a good initialization in the nonconvex problem.

\paragraph{Primal feasibility}
In approximate dual methods, primal feasibility is usually not assured. Here, we just use the maximum eigenvector of $\mA^*(y^*)$ to reconstruct the image, but first performing some projection to ensure primal feasibility may lead to better answers. 

\paragraph{Other spectral approximation methods} Here, we experiment with a spectral approximation method unique to the phase retrieval problem, where in the limit, $\mA^*(\mA(x)) \approx x$. We have not compared this to other spectral approximation methods, like sketching, subsampling, sparsification, etc.

\bibliographystyle{alpha}
\bibliography{refs}

\end{document}